\documentclass[12pt]{amsart}

\usepackage{graphicx} 
\usepackage{setspace} 
\usepackage{amsfonts} 
\usepackage{amsmath} 
\usepackage{amssymb}
\usepackage[margin=1in]{geometry}
\usepackage{amsthm}
\usepackage{bbm}
\usepackage{amssymb,latexsym,amsmath}
\newtheorem{theorem}{Theorem}[section]

\newtheorem{remark}[theorem]{Remark}

\numberwithin{equation}{section} 
\begin{document}
%
\title{Phase Transition in the Recovery of Rank One Matrices Corrupted by Gaussian Noise}

\author{Enrico Au-Yeung} 
\address{Department of Mathematical Sciences \\
DePaul University\\
Chicago \ IL, 60614 }
\email{eauyeun1@depaul.edu}

\author{Greg Zanotti}
\address{Department of Mathematical Sciences \\
DePaul University\\
Chicago \ IL, 60614 }
\email{gregzanotti@gmail.com}


\maketitle
\date{}

\doublespacing

\begin{abstract}
 In datasets where the number of parameters is fixed and the number of samples  is large, principal component analysis (PCA) is a powerful dimension reduction tool.  However, in many contemporary datasets, when the number of parameters is comparable to the sample size, PCA can be misleading. A closely related problem is the following: is it possible to recover a rank-one matrix in the presence of a large amount of noise?  In both situations, there is a phase transition in the eigen-structure of the matrix.
 
\end{abstract}


\section{Introduction}\label{intro}

The problem of low-rank matrix recovery has received a lot of attention in the signals processing community over the last 10 years.  One reason is the close connection to the matrix completion problem: 
Given a random subset of the entries of a low-rank matrix, the task is to recover all the entries of this matrix.  According
to the foundational paper by Candes and Recht \cite{CandesRecht2009}, this task can be accomplished by solving a convex optimization problem. 
    The matrix completion problem arises in machine learning, image processing, and the Netflix-Prize problem.  The practical nature of this problem has motivated many researchers to investigate efficient methods to solve this optimization problem.
See, e.g., \cite{CaiCandesShen2010}, \cite{GoldfarbMa2011}, \cite{KeshavanMontanariOh2010},   \cite{Vandenberghe2009}, \cite{MohanFazel2012}, \cite{FornasierRauhutWard2011},
\cite{ChenHeYuan2012}, \cite{WenYinZhang2012}, \cite{Chatterjee2015}.  This list is by no means exhaustive. 

In this paper, we take a different direction.  We are interested in the following type of situation: Is it possible to recover a rank-one matrix in the presence of a large amount of noise?  A useful data model to keep in mind is the following: $X = \lambda \bf{x} \bf{x}^{T} + G$, where $\bf{x} \in \mathbb{R}^n$.  The data matrix $X$ represents our observations, and the Gaussian matrix $G$ represents the noise structure.  The challenge is to recover the principal vector $\bf{x}$ and the value $\lambda$ from the data matrix $X$.  We are especially interested in  the asymptotic behaviour of the largest eigenvalue and the leading eigenvector  of the data matrix $X$, (as $n \rightarrow \infty$), when the operator norm of $G$ is not negligible compared to the operator norm of $X$.  We observe a phase transition phenomenon.

We are also interested in the behaviour of the leading singular vector when the data matrix is a large rectangular matrix, where the number of rows is proportional to the number of columns.  
Principal component analysis (PCA)  is a versatile tool in dimensionality reduction. PCA projects the data onto the principal subspace 
spanned by the leading eigenvectors of the sample covariance matrix.  In theory, these eigenvectors can capture most of the variance in the data.  This enables the dimension of the feature space to be reduced, while retaining most of the
information.  In the contemporary setting, a collection of high-dimensional data can be treated as a low-rank signal with additional noise
structure.  If the samples of data are organized into a data matrix, then PCA can be used to recover the low-rank signal.  It performs well
when the number of features, $p$, is small, and the number of samples $n$ is large (\cite{Anderson1963}).  However, in biomedical studies, the number of features $p$ is often comparable to the sample size $n$.  In the biomedical setting, the features are measurements on the expression levels of thousands of genes, and $n$ is the thousands of individuals.

\subsection{Setting and Motivation} 
Suppose we have a collection of independently and identically distributed random vectors, $x_1, x_2, x_3, \ldots, x_n$  from a $p$-dimensional real Gaussian distribution with mean zero and covariance $\Sigma =\mbox{diag}(\lambda_1, \lambda_2, \ldots, \lambda_M, 1, 1, \ldots, 1)$, where $\lambda_1 \geq \lambda_2 \geq \ldots \geq \lambda_M > 1$.  Let $X$ be the $p \times n$ matrix with column vectors $x_1, x_2, \ldots, x_n \in \mathbb{R}^p$. Assume that $0 < c < 1$ and $\frac{p}{n} = c$.  Let $S = \frac{1}{n} X X^{T}$ be the sample covariance matrix.  A data scientist wants to know how the largest eigenvalues of the matrix $S$ behave as $n \rightarrow \infty$.  Let us consider a specific scenario, when $p = 500$ and $n = 2000$,  is the sample largest eigenvalue $\widehat{\lambda}_1$ of the matrix $S$ a good estimator of the true eigenvalue $\lambda_1$? That depends on the true value of the largest eigenvalue $\lambda_1$. 

The following is a simplified version of a theorem of Baik and Silverstein (see \cite{BaikSilverstein2006}, \cite{Paul2007}): 
\begin{theorem}

  Let $\widehat{\lambda}_1$ be the largest eigenvalue of $S$.\\
(1).  Suppose $\lambda_1 \leq 1 + \sqrt{c}$ and $\frac{p}{n} \rightarrow c$ as $n \rightarrow \infty$.  Then,  we have \[ \widehat{\lambda}_1 \rightarrow (1 + \sqrt{c})^2 \quad \mbox{ as } n \rightarrow \infty. \]
(2).  Suppose $\lambda_1 > 1 + \sqrt{c}$ and $\frac{p}{n} \rightarrow c$ as $n \rightarrow \infty$.  Then, we have \[ \widehat{\lambda}_1 \rightarrow \lambda_1 (1 + \frac{c}{\lambda_1 - 1}) \mbox{ as } n \rightarrow \infty. \]

\end{theorem}
What these authors observe is that there is a phase transition in the eigen-structure of a matrix when both the rows and columns are large, i.e. when $\frac{p}{n} \rightarrow c, 0 < c < 1$ and $n \rightarrow \infty$. \\  The phase transition phenomenon can be quite complicated and this has been analyzed in the seminal paper \cite{BBP2005}.  For other variations on this theme, see, e.g. \cite{Peche2006}, \cite{FeralPeche2007}, and \cite{PizzoRenfrewSoshnikov2013}.\\

  In the data model, $X = \lambda \bf{x_1} \bf{x_1}^{T} + G$, suppose we have some additional information about the principal vector $\bf{x}$, how can we use that information to recover the vector?   We consider the case when each entry of the vector is bounded between $0$ and a known constant $\tau$.  To be precise, we can take the specific value, $\tau = 0.2$. Thus, we know that the vector $\bf{x_1}$ lies in a box.   
  
   Numerical experiments show that, when $\lambda = 4$, the leading eigenvector $\bf{v_1}$ of the data matrix $X$ is  not a good approximation to the desired vector $\bf{x_1}$.  In fact, the relative error between $\bf{v_1}$ and $\bf{x_1}$ often exceeds a hundred percent.  
   
   The purpose of this paper is to address both the theoretical and practical aspect of this problem.  On the theory side, we observe a phase transition in the largest eigenvalue.  Moreover, the result shows that, depending on the true value of $\lambda_1$, the leading eigenvector of the data matrix $X$ can be nearly orthogonal to the true vector $\bf{x_1}$.  This means that some caution is warranted: when $n \rightarrow \infty$, using principal component analysis as an attempt to retrieve the vector $\bf{x_1}$ can give a misleading result. 
   In place of a proof to the theorem, we provide a heuristic explanation that gives the main insight to the theorem.\\
   
   
    On the practical side, 
  we develop an iterative algorithm to recover the vector $\bf{x_1}$ from the matrix $X$.  We view this as a box-constrained optimization problem to find a vector ${\bf{x}}$, 
  where the unknown variable satisfies the constraint, $0 \leq \|\bf{x}\|_{\infty} \leq \tau$.
   Compared to the leading eigenvector ${\bf{v_1}}$ of the matrix $X$, our algorithm yields a vector that is significantly closer to the desired vector ${\bf{x_1}}$. \\
But first, we need to set some
Notations: \\
$X$ is a symmetric random matrix,  $X \in \mathbb{R}^{n \times n}$ \\
${\bf{x}_1}$ is a fixed (non-random) vector, $\| {\bf{x_1}} \|_2 = 1,$ and ${\bf{x_1}} \in \mathbb{R}^n$ \\
$G$ is a Gaussian symmetric matrix,  $G \in \mathbb{R}^{n \times n}$ and $G = G^{T}$, where \\
 $G(i,j)$ are independent, normally distributed with mean 0 and variance $\frac{1}{n}$ for $i < j$, and
$G(i,i)$ is normally distributed with mean 0 and variance $\frac{2}{n}$ \\

The following theorem is the phase transition phenomenon (for symmetric matrices).\\
{\bf{Important Note Added:}} After the initial preparation of an earlier version of this manuscript, we learned that this is a version of a theorem of Florent Benaych-Georges and Raj Rao Nadakuditi, see \cite{FlorentRajRao2011}. We encourage the reader to consult this beautifully written paper.  We are grateful to the authors of that paper for bringing this to our attention.

\begin{theorem}\label{thm:symmetric_Gaussian}

Let $X = \lambda {\bf{x_1}} {\bf{x_1}}^{T} + G$, where G is Gaussian symmetric matrix.  Pick $\tau = 0.2$.   \\
Suppose ${\bf{x_1}}$ is a fixed vector of length 1, and $0 \leq {\bf{x_1}}(j) \leq \tau$, for $1 \leq j \leq n$. \\
Let $\widehat{\lambda}_1$ be the largest eigenvalue of the matrix $X$.
Let ${\bf{v_1}}$ be the leading eigenvector of the matrix $X$, i.e.  ${\bf{v_1}}$ is the eigenvector that corresponds to $\widehat{\lambda}$.

Then, if $\lambda > 1$, we have 
\[ \lim_{n  \rightarrow \infty} | \langle {\bf{v_1}}, {\bf{x_1}} \rangle | = \sqrt{c}, \]
where $c = 1 - \frac{1}{\lambda^2}$.  Otherwise, if $\lambda \leq 1$, we have
\[ \lim_{n  \rightarrow \infty} | \langle {\bf{v_1}}, {\bf{x_1}} \rangle | = 0. \]
For the largest eigenvalue of the matrix $X$, the following phase transition occurs.
If $\lambda \geq 1$,  we have \begin{equation}\label{eqn:lambda_plus}\widehat{\lambda}_1 \rightarrow \lambda + \frac{1}{\lambda}\end{equation} as $n \rightarrow \infty$. Otherwise, if $\lambda \leq 1$, we have $\widehat{\lambda}_1 \rightarrow 2$ as $n \rightarrow \infty$.\\

\end{theorem}
\begin{remark}Interestingly, in their theorem \cite{FlorentRajRao2011}, they do not have the hypothesis that the vector ${\bf{x_1}}$ satisfied the constraint $0 \leq {\bf{x_1}}(j) \leq \tau$, for $1 \leq j \leq n$, and the conclusion of their theorem remains the same.  We include this additional condition in the statement of the theorem, since we explicitly use this condition in our numerical optimization algorithm.\end{remark}

\section{Background for Wigner matrices}

The symmetric Gaussian random matrix in Theorem \ref{thm:symmetric_Gaussian}  is an example of a Wigner matrix. 
We summarize here some background and standard facts regarding the Wigner Semicircular Law for symmetric random matrices.    Given any probability measure on the real line, the Stieltjes transform  is defined by
\[ S_\mu(z) = \int_{\mathbb{R}} \frac{d\mu(t)}{z - t},\]
where $z$ is any complex number  in the upper half of the complex plane.  For any $n \times n$ symmetric matrix $M_n$, we can work with the normalized matrix   $\frac{1}{\sqrt{n}}M_n$ and form its empirical spectral distribution (ESD),
\[ \mu_{\frac{1}{\sqrt{n}}M_{n}}(x) = \frac{1}{n} \sum_{j=1}^{n} \delta\left(x - \frac{\lambda_j(M_n)}{\sqrt{n}}\right)\] of $M_n$, where $\lambda_j(M)$ are the  eigenvalues of $M_n$.  The ESD is a probability measure, also known as the spectral measure for the matrix.  For the square matrix $M_n$ with spectral measure $\mu(x) = \mu_{M_n}(x)$, we can define its corresponding Stieltjes transform.  We have the following useful identity,
\[S_n(z) = S_{\mu_{\frac{1}{\sqrt{n}}M_{n}}}(z)  = \frac{1}{n} Tr\left[\left( \frac{1}{\sqrt{n}}M_{n} - z I_{n}\right)^{-1}\right]\]
where $Tr$ denotes the trace of a matrix, and $I_n$ is the $n \times n$ identity matrix.  We define the semicircular distribution $\mu_{sc}(x) = \frac{1}{2 \pi} \sqrt{4 - x^2}$.  The Wigner semicircular law states that the sequence of ESDs $\mu_{\frac{1}{\sqrt{n}}M_{n}}(x)$ converges almost surely to $\mu_{sc}(x)$.  The Stieltjes transform for the spectral measure $\mu_{sc}$ is 
\[ S_{\mu_{sc}}(z) = \int_{\mathbb{R}} \frac{d \mu_{sc}(x)}{x - z} = \frac{-z + \sqrt{z^2 - 4}}{2}.\]

\section{Main insight for Theorem \ref{thm:symmetric_Gaussian}}

We now give the heuristic explanation for the quantity $\lambda + \frac{1}{\lambda}$  that appears in equation (\ref{eqn:lambda_plus})   in the phase transition phenomenon of Theorem \ref{thm:symmetric_Gaussian}. \\

{\bf{Note Added:}} After the initial preparation of an earlier version of this manuscript, we learned that  there is a similar discussion  in  \cite{FlorentRajRao2011}, and is accompanied by a rigorous proof.  We are grateful to the authors of that paper for bringing this to our attention.  \\

Recall that $X = \lambda {\bf{x_1}} {\bf{x_1}}^{T} + G$, where G is Gaussian symmetric matrix. Since $G$ is symmetric, we can write $G = U^T D U$, where $U$ is an orthogonal matrix and $D = \mbox{diag}(\lambda_1, \lambda_2, \ldots, \lambda_n)$ is a diagonal matrix.  Instead of the matrix $X$, we can consider the matrix $U X U^{T} = D + \lambda U {\bf{x_1}} {\bf{x_1}}^{T} U^T$, i.e.  a diagonal matrix $D$ plus a rank one positive matrix $P \equiv \lambda U {\bf{x_1}} {\bf{x_1}}^{T} U^T$.  The random orthogonal matrix $U$ rotates the fixed vector $\bf{x}$ of length one to a random vector $\bf{u}$ of length one.  The intuition is that when $n$ is large, then with high probability, the vector $\bf{u}$ is uniformly distributed on the unit sphere $\{x \in \mathbb{R}^n \colon \| x\|_2 = 1\}.$  Hence, each entry $u(k)$ of the unit-length vector $\bf{u}$ is approximately equal to the square root of $1/n$.  Fix $z$ and suppose the matrix $(D - z I_n)$ is invertible.  Then, we have the relation,
\begin{equation} 
\det(zI_n - (D + P)) = \det( zI_n - D) \cdot \det(I_n - (z I_{n} - D)^{-1} P). 
\end{equation}
Consider the matrix $M \equiv (z I_n - D)^{-1} P.$  Then, 1 is an eigenvalue of the matrix $M$ if and only if $z$ is not an eigenvalue of $D$ and $z$ is an eigevalue of $D + P$.  Since the matrix $M = (z I_n - D)^{-1} \lambda \bf{u}\bf{u}^{T}$ has rank one, so the trace of $M$ is equal to the only nonzero eigenvalue of $M$.  On the other hand, we have
\[Tr(M) = \lambda \sum_{k=1}^{n} \frac{| u(k)|^2}{z  - \lambda_{k}}.\]   
This implies that $z$ is not an eigenvalue of $D$ and $z$ is an eigenvalue of $D + P$ if and only if 
\begin{equation}\label{eqn:spectral_D_plus_P}{\lambda} \sum_{k=1}^{n} \frac{ | u(k) |^2 }{z - \lambda_k} = 1. \end{equation}
Here, $u(k)$ are the entries of the vector $\bf{u}$.
The left hand side of (\ref{eqn:spectral_D_plus_P}) is $\lambda S_{\mu_n}(z)$, where $\mu_n$ represents a weighted spectral measure associated to the diagonal matrix $D$, \[ \mu_n(x) = \sum_{k=1}^{n} |u(k)|^2 \cdot \delta(x - \lambda_k).\] 
Recall that when $n$ is large, the square of each entry $u(k)$ of the vector $\bf{u}$ is about $1/n$, with high probability.  Thus, we replace the previous relation (\ref{eqn:spectral_D_plus_P}) with
\begin{equation}\label{eqn:spectral_Wigner} \frac{1}{n} \sum_{k=1}^{n} \frac{1}{z - \lambda_k} = \frac{1}{\lambda}. \end{equation}
But the left hand side of equation (\ref{eqn:spectral_Wigner}) converges to the Stieltjes transform of the semicircular distribution, so we conclude that
\[S_{\mu_{sc}}(z) = \frac{1}{\lambda}.\]
Inverting the Stieltjes transform, we have $z = S_{\mu_{sc}}^{-1}(\frac{1}{\lambda})$. By direct calculation, we can verify that
\[ S_{\mu_{sc}}^{-1}(\frac{1}{\lambda}) = \lambda + \frac{1}{\lambda}. \]
Finally, since $z$ is an eigenvalue of $D + P$, we have shown that $\lambda + \frac{1}{\lambda}$ is indeed an eigenvalue of $D + P$.  This completes our heuristic explanation for the appearance of $\lambda + \frac{1}{\lambda}$ in the phase transition of eigenvalues for symmetric matrices.

\section{Optimization Algorithm}
Our optimization algorithm is based on gradient descent, but includes additional steps to satisfy our problem's constraints. 
The input is the symmetric matrix \textbf{X}, which is the observed data with Gaussian noise structure.
We use an iterative algorithm to estimate the unknown vector $\bf{x}$, under the constraint (where $\tau$ is an input parameter):
\begin{equation}
\begin{aligned}
\quad & 0 \leq \mathbf{x}(i) \leq \tau \quad \forall i \in \{1, 2, ..., N\} \\
\quad & ||\mathbf{x}||_2 = 1.
\end{aligned}
\end{equation}
As the additive Gaussian noise increases the magnitude of the eigenvalues of the observed matrix \textbf{X}, we perform gradient descent on an estimate of the true gradient (i.e.~for $\mathbf{x_1x_1^T}$ instead of \textbf{X}) by penalizing the magnitude of the eigenvalues of our recovered matrix: that is, we penalize the trace of $\mathbf{x}\mathbf{x}^T$ by adding as a penalty term the $L_2$ norm of \textbf{x}. This penalty helps minimize the nuclear norm of $\mathbf{x}\mathbf{x}^T$. Our update equation is as follows:
\begin{equation}
\begin{aligned}
\mathbf{x_{k+1}} = \mathbf{x_k} - \alpha [\mathbf{x_k^T(x_k x_k^T - X)} + \gamma (\mathbf{x_k^T}\mathbf{x_k})\mathbf{1^T}]^{\mathbf{T}}\\
\end{aligned}
\end{equation}
Above, $\alpha$ is the usual step size, which we experimentally observe to work best when set in the range $[10^{-3}, 10^{-1}]$, with slightly better recovery results toward the lower end of the range. The second parameter $\gamma$ is the regularization parameter for the $L_2$ norm of \textbf{x}, which we set to $10^{-1}$ and do not change. Our gradient descent continues until the following termination condition is met, which is usually satisfied within roughly 50 iterations when $\alpha = 0.1$:
\begin{equation}
\begin{aligned}
\frac{||\mathbf{x_{k+1}} - \mathbf{x_k}||_2}{||\mathbf{x_{k+1}}||_2} \leq 10^{-5}.
\end{aligned}
\end{equation}
We initialize $\mathbf{x}$ by generating a length $N$ vector of uniform random numbers in $[0, \tau]$ and then dividing it by its $L_2$ norm.

This gradient descent procedure, however, does not account for the main constraint in our optimization problem: the box constraint. To satisfy this constraint, after gradient descent reaches the termination condition above, we apply a projection step that mitigates the effect of the additional additive noise in the off-diagonal entries of the observed matrix \textbf{X}. First, we divide $\mathbf{x}$ by its $L_2$ norm. Then, we project $\mathbf{x}$ onto the box $[0, \tau]^n$ by setting each $\mathbf{x}(i) = \text{min}(\text{max}(\mathbf{x}(i), 0), \tau)$. Finally, we again divide $\mathbf{x}$ by its $L_2$ norm.

In our experiments, we initialize the true vector $\mathbf{x_1}$ by setting a block of 2\% of the entries to $1 / \sqrt{2(10^{-2})n}$ and dividing it by its $L_2$ norm. We keep $\alpha=10^{-1}$. The observed data \textbf{X} is  the rank-one matrix $\lambda \bf{x_1x_1^{T}}$,  plus a Gaussian random matrix $G$, as described in the previous section, with $\lambda=4$. We define the relative error as:
\begin{equation}
\begin{aligned}
\text{E}(\mathbf{x}) = 100 \cdot \frac{||\mathbf{x_{1}} - \mathbf{x}||_2}{||\mathbf{x_1}||_2}
\end{aligned}
\end{equation}
where \textbf{x} is our recovered vector. Regardless of our selection of $\alpha$ in the range above, the standard deviations of the relative error at each $n$ are consistently below $2\%$ for $n \geq 500$. For sizes of $\mathbf{x_1}$ ranging from $500$ to $5000$, we observe that the average relative error for the recovered vector using our optimization is substantially lower in comparison to that of using the leading eigenvector.  Average relative error is computed over 200 trials (i.e.~draws of $G$ and optimization procedures) per $n$. The  results of our experiment are displayed in Table \ref{tbl:experiments}.

\newcommand\tstrut{\rule{0pt}{2.5ex}}
\newcommand{\avgerr}{\text{Mean}(E)}
\newcommand{\sderr}{\text{SD}(E)}

\begin{table}[]
	\centering
	\caption{Relative error averaged over 200 trials. $\avgerr$ denotes the mean relative error. ``Opt'' is our optimization procedure and ``Eig'' is the top eigenvector procedure.}
	\begin{tabular}{cccccccc}
	\hline \tstrut
	$n$  & Opt $\avgerr$ & Eig $\avgerr$            \\ \hline
	500  & 15.4\% 		  & 113.5\%		   \tstrut  \\ \hline
	1000  & 14.3\% 		  & 107.4\%		   \tstrut  \\ \hline
	2500  & 12.5\% 		  & 123.9\%		   \tstrut  \\ \hline
	5000  & 10.6\% 		  & 111.8\%		   \tstrut  \\ \hline
	\end{tabular}
	\label{tbl:experiments}
\end{table}


%



\bibliographystyle{plain}
\bibliography{PhaseTransitionSAMPTARefs}

\begin{thebibliography}{10}

\bibitem{Anderson1963}
T.~W. Anderson.
\newblock Asymptotic theory for principal component analysis.
\newblock {\em Ann. Math. Statist.}, 34:122--148, 1963.

\bibitem{BBP2005}
Jinho Baik, G\'erard Ben~Arous, and Sandrine P\'ech\'e.
\newblock Phase transition of the largest eigenvalue for nonnull complex sample
  covariance matrices.
\newblock {\em Ann. Probab.}, 33(5):1643--1697, 2005.

\bibitem{BaikSilverstein2006}
Jinho Baik and Jack Silverstein.
\newblock Eigenvalues of large sample covariance matrices of spiked population
  models.
\newblock {\em J. Multivariate Anal.}, 97(6):1382--1408, 2006.

\bibitem{FlorentRajRao2011}
Florent Benaych-Georges and Raj~Rao Nadakuditi.
\newblock The eigenvalues and eigenvectors of finite, low rank perturbations of
  large random matrices.
\newblock {\em Adv. in Math.}, 227:494--521, 2011.

\bibitem{CaiCandesShen2010}
J.F. Cai, Emmanuel~J. Candes, and Z.W. Shen.
\newblock A singular value thresholding algorithm for matrix completion.
\newblock {\em SIAM J. Optim.}, 20:1956--1982, 2010.

\bibitem{CandesRecht2009}
Emmanuel~J. Candes and Benjamin Recht.
\newblock Exact matrix completion via convex optimization.
\newblock {\em Found. Comput. Math.}, 9(6):717--772, 2009.

\bibitem{Chatterjee2015}
Sourav Chatterjee.
\newblock Matrix estimation by universal singular value thresholding.
\newblock {\em Ann. Statist.}, 43(1):177--214, 2015.

\bibitem{ChenHeYuan2012}
Caihua Chen, Bingsheng He, and Xiaoming Yuan.
\newblock Matrix completion via an alternating direction method.
\newblock {\em IMA J. Numer. Anal.}, 32(1):227--245, 2012.

\bibitem{FeralPeche2007}
Delphine F\'eral and Sandrine P\'ech\'e.
\newblock The largest eigenvalue of rank one deformation of large wigner
  matrices.
\newblock {\em Comm. Math. Phys.}, 272(1):185--228, 2007.

\bibitem{FornasierRauhutWard2011}
Massimo Fornasier, Holger Rauhut, and Rachel Ward.
\newblock Low-rank matrix recovery via iteratively reweighted least squares
  minimization.
\newblock {\em SIAM J. Optim.}, 21(4):1614--1640, 2011.

\bibitem{GoldfarbMa2011}
Donald Goldfarb and Shiqian Ma.
\newblock Convergence of fixed-point continuation algorithms for matrix rank
  minimization.
\newblock {\em Found. Comput. Math.}, 11(2):183--210, 2011.

\bibitem{KeshavanMontanariOh2010}
Raghunandan~H. Keshavan, Andrea Montanari, and Sewoong Oh.
\newblock Matrix completion from a few entries.
\newblock {\em IEEE. Trans. Inf. Theory}, 56(6):2980--2998, 2010.

\bibitem{Vandenberghe2009}
Z.~Liu and L.~Vandenberghe.
\newblock Interior-point method for nuclear norm approximation with application
  to system identification.
\newblock {\em SIAM. J. Matrix Anal. Appl.}, 31:1235--1256, 2009.

\bibitem{MohanFazel2012}
Karthik Mohan and Maryam Fazel.
\newblock Iterative reweighted algorithms for matrix rank minimization.
\newblock {\em J. Mach. Learn. Res.}, 13:3441--3473, 2012.

\bibitem{Paul2007}
Debashis Paul.
\newblock Asymptotics of sample eigenstructure for a large dimensional spiked
  covariance model.
\newblock {\em Statist. Sinica}, 17(4):1617--1642, 2007.

\bibitem{Peche2006}
Sandrine P\'ech\'e.
\newblock The largest eigenvalue of small rank perturbations of hermitian
  random matrices.
\newblock {\em Probab. Theory Related Fields}, 134(1):127--173, 2006.

\bibitem{PizzoRenfrewSoshnikov2013}
Alessandro Pizzo, David Renfrew, and Alexander Soshnikov.
\newblock On finite rank deformations of wigner matrices.
\newblock {\em Ann. Inst. Henri Poincare Probab. Stat.}, 49(1):64--94, 2013.

\bibitem{WenYinZhang2012}
Zaiwen Wen, Wotao Yin, and Yin Zhang.
\newblock Solving a low-rank factorization model for matrix completion by a
  nonlinear successive over-relaxation algorithm.
\newblock {\em Math. Program. Comput.}, 4(4):333--361, 2012.

\end{thebibliography}
%
%
%


\end{document}